\input amstex
\magnification=\magstep1 
\baselineskip=13pt
\documentstyle{amsppt}
\vsize=8.7truein \CenteredTagsOnSplits \NoRunningHeads
 \def\per{\operatorname{per}}
 \def\PP{\bold{Pr}\thinspace}
 \def\EE{\bold{E}\thinspace}

 \topmatter

\title  Computing the partition function for perfect matchings in a hypergraph \endtitle
\author Alexander Barvinok and Alex Samorodnitsky \endauthor
\address Department of Mathematics, University of Michigan, Ann Arbor,
MI 48109-1043, USA \endaddress
\email barvinok$\@$umich.edu \endemail
\address Department of Computer Science, Hebrew University of Jerusalem, Givat Ram Campus,
91904, Israel \endaddress
\email salex$\@$cs.huji.ac.il \endemail
\date September 2011 \enddate
\thanks  The research of the first author was partially supported by NSF Grant DMS 0856640.
The research of the second author was partially supported by ISF grant 039-7165.
The research of the authors was also partially supported 
by a United States - Israel BSF grant 2006377.
\endthanks 

\abstract Given non-negative weights $w_S$ on the $k$-subsets $S$ of a $km$-element set $V$, we consider 
the sum of the products $w_{S_1} \cdots w_{S_m}$ over all partitions 
$V=S_1 \cup \ldots \cup S_m$ into pairwise disjoint $k$-subsets $S_i$. When the weights $w_S$ are positive and within a constant factor, fixed in advance, of each other, we present a simple polynomial time algorithm to approximate the sum within a polynomial in $m$ factor. In the process, we obtain higher-dimensional versions of the van der Waerden and 
Bregman-Minc bounds for permanents. We also discuss applications to counting of perfect and nearly perfect matchings
in hypergraphs.
\endabstract
\endtopmatter

\document

\head 1. Introduction and main results \endhead

Let us fix an integer $k>1$. A collection $H \subset {V \choose k}$ of $k$-subsets of a finite set $V$ is called
a {\it $k$-uniform hypergraph} with {\it vertex} set $V$, while sets $S \in H$ are called {\it edges} of $H$.
In particular, a uniform $2$-hypergraph is an ordinary undirected 
graph on $V$ without loops or multiple edges. A set $\left\{S_1, \ldots, S_m \right\}$ of pairwise vertex disjoint 
edges of $H$ such that $V=S_1 \cup \ldots \cup S_m$ is called a {\it perfect matching} of hypergraph $H$.
More generally, a {\it matching} of {\it size} $n$ is a collection of $n$ pairwise disjoint edges of $H$.

If a perfect matching exists then the number $|V|$ of vertices of $V$ is divisible by $k$, so we have 
$|V|=km$ for some integer $m$. The hypergraph consisting of all $k$-subsets of $V$ is called the
{\it complete $k$-uniform hypergraph} with vertex set $V$. We denote it by ${V \choose k}$. 
A hypergraph is called a {\it complete $k$-partite hypergraph} if the set $V$ of vertices is a union 
$V=V_1 \cup \ldots \cup V_k$ of pairwise disjoint sets $V_i$, called {\it parts}, such that $\left|V_1\right| =\ldots =\left|V_k\right|=m$
and the edges of the hypergraph are the subsets $S \subset V$ containing exactly one vertex in each part:
$\left|S \cap V_1\right| =\ldots =\left|S \cap V_k \right|=1$. We denote such a hypergraph by $V_1 \times \ldots \times V_k$.

We introduce the main object of the paper.

\definition{(1.1) Partition function} Let $H$ be a $k$-uniform hypergraph with the set $V$ of vertices
such that $|V|=km$ for some positive integer $m$. Suppose that to every edge 
$S \in H$ a non-negative real number $w_S$ is assigned. Such an assignment $W=\left\{w_S \right\}$ we call
a {\it weight} on $H$. We say that $W$ is {\it positive} if $w_S>0$ for all $S \in H$. The polynomial
$$P_H(W) =\sum  w_{S_1} \cdots w_{S_m}, $$
where the sum is taken over all perfect matchings $\left\{ S_1, \ldots, S_m \right\}$ of $H$, is called 
the {\it partition function} of perfect matchings in hypergraph $H$. Sometimes we write just $P(W)$ if the choice of 
the hypergraph $H$ is clear from the context.

We note that we can obtain the partition function $P_H(W)$ of an arbitrary $k$-uniform hypergraph $H \subset {V \choose k}$ by 
specializing $w_S=0$ for $S \notin H$ in the partition function of the complete $k$-uniform 
hypergraph ${V \choose k}$. 

The partition function of ${V \choose 2}$ with $|V|=2m$ is known as the {\it hafnian} of the $2m \times 2m$ symmetric matrix 
$A=\left(a_{ij}\right)$, where $a_{ij}$ is the weight of the edge consisting of the $i$-th and $j$-th vertices of $V$ (diagonal elements of $A$ can be chosen arbitrarily), see, for example, Section 8.2 of \cite{Mi78}. If $V_1 \times V_2$ is a complete 
bipartite graph with $\left|V_1\right| =\left| V_2 \right| =m$ then the corresponding partition function is the {\it permanent} of 
the $m \times m$ matrix $B=\left(b_{ij}\right)$, where $b_{ij}$ is the weight of the edge consisting of the $i$-th vertex
of $V_1$ and $j$-th vertex of $V_2$.
The partition function of the complete $k$-partite hypergraph
gives rise to a version of the permanent of a $k$-dimensional tensor, see, for example, \cite{D87b}.
\enddefinition

In this paper, we address the problem of computing or approximating $P_H(W)$ efficiently. First, we define certain classes of  weights $W$.

\definition{(1.2) Balanced and $k$-stochastic weights} We say that a positive weight $W=\left\{w_S \right\}$ on 
a $k$-uniform hypergraph is {\it $\alpha$-balanced} for 
some $\alpha \geq 1$ if 
$${w_{S_1} \over w_{S_2}} \ \leq \ \alpha \quad \text{for all} \quad S_1, S_2 \in H.$$
Note that an $\alpha$-balanced weight is also $\beta$-balanced for any $\beta > \alpha$.

Weight $Z=\left\{z_S\right\}$ is called {\it $k$-stochastic}, if 
$$\sum \Sb S \in H\\  S \ni v \endSb z_S=1 \quad \text{for all} \quad v \in V.$$
In words: for every vertex, the sum of the weights of the edges containing the vertex is 1.
\enddefinition

Now we are ready to state our first main result.

\proclaim{(1.3) Theorem} Let us fix an integer $k>1$ and a real $\alpha \geq 1$. Then there exists
a real $\gamma=\gamma(k, \alpha) >0$ such that if $H$ is a complete $k$-uniform hypergraph 
or a complete $k$-partite hypergraph with $km$ vertices and $Z$ is 
a $k$-stochastic $\alpha$-balanced weight on $H$ then 
$$m^{-\gamma} e^{-m(k-1)} \ \leq \ P_H(Z) \ \leq \ m^{\gamma} e^{-m(k-1)}$$
provided $m>1$.
\endproclaim
In other words, for fixed $k$ and $\alpha$, the value of the partition function for a $k$-stochastic $\alpha$-balanced weight 
on a complete $k$-uniform hypergraph or a complete $k$-partite hypergraph can vary only within a polynomial in $m$ range. 

More precisely, we prove that under conditions of Theorem 1.3 and assuming, additionally, that $\alpha^{k+1}>2$, we have 
$$\epsilon_1 m^{-\gamma_1} e^{-m(k-1)}  \ \leq \ P_H(Z) \ \leq \ \epsilon_2 m^{\gamma_2} e^{-m(k-1)},$$
where 
$$\aligned &\gamma_1 = \alpha^{3(k+1)}(k^2+k)^2 +(k-1)^2 \quad \text{and} \quad  \gamma_2 = {k^2 \alpha^{k+1} \over 2},\\
&\epsilon_1 =\alpha^{-(k+1)l} l^{l} {kl \choose k}^{1-l} \quad \text{and} \quad  \epsilon_2=\alpha^{(k+1)l} l^{l -kl+k} \\
&\qquad \qquad \text{for} \\
& l=\lceil \alpha^{2(k+1)} k^2 \rceil +1.\endaligned \tag1.3.1$$
\subhead (1.4) Comparison with permanents \endsubhead
The van der Waerden conjecture on permanents proved by Falikman \cite{Fa81} and Egorychev \cite{Eg81}, 
see also \cite{Gu08} for important new developments, 
asserts that if $A=\left(a_{ij}\right)$ is an $m \times m$ doubly stochastic matrix, that is, a non-negative matrix with all row 
and column sums equal 1, then
$$\per A \ \geq \ {m! \over m^m} =\sqrt{2 \pi m} e^{-m} \left(1 +O\left({1 \over m}\right)\right).$$
A conjecture by Minc proved by Bregman \cite{Br73}, see also \cite{Sc78} for a simpler proof, asserts that
if $B=\left(b_{ij}\right)$ is an $m \times m$ matrix with $b_{ij} \in \{0, 1\}$ for all $i,j$ then 
$$\per B \ \leq \ \prod_{i=1}^m \left(r_i ! \right)^{1/r_i},$$
where $r_i$ is the $i$-th row sum of $B$. From this inequality one can deduce that if $A$ is an $m \times m$ non-negative 
matrix with all row sums equal 1 and all the entries not exceeding $\alpha/m$ for some 
$\alpha \geq 1$ then 
$$\per A \ \leq \ m^{\gamma} e^{-m}$$
for some $\gamma=\gamma(\alpha)>0$ and all $m>1$ (one can choose any $\gamma > \alpha/2$ 
if $m$ is sufficiently large), see \cite{So03}.
Thus the van der Waerden and Bregman-Minc inequalities together imply that $\per A =e^{-m} m^{O(1)}$ for any 
$m \times m$ doubly stochastic matrix $A$ whose entries are within a factor of $O(1)$ of each other. Theorem 1.3
presents an extension of this interesting fact to non-bipartite graphs for $k=2$ and to hypergraphs 
for $k >2$. A stronger statement that $\per A =e^{-m} m^{O(1)}$ for an $m \times m$ doubly stochastic matrix whose 
maximum entry is $O\left(m^{-1}\right)$ fails to extend to non-bipartite graphs for $k=2$ or to $k$-partite 
hypergraphs for $k>2$ as the following two examples readily show.

Let $k=2$ and let $H$ be a graph on a set $V$ of $n=4r+2$ vertices, which consists of two vertex-disjoint copies of 
the complete graph on 
$2r+1$ vertices. Let us define a weight $Z=\left\{z_S\right\}$ on ${V \choose 2}$ by letting $z_S=(2r)^{-1}$ if $S$ is an edge of $H$ and $z_S=0$ otherwise. Then $Z$ is 2-stochastic weight on ${V \choose 2}$ and $P(Z)=0$. That is, the hafnian 
of an $n \times n$ symmetric doubly stochastic matrix can be zero even when the maximum entry of the matrix is 
$O\left(n^{-1}\right)$.

Let $k=3$, let $m=4r+2$ and let us identify each set $V_1, V_2$ and $V_3$ with a copy of the 
integer interval $\{1, 2, \ldots, m\}$. Let us define a 
weight $Z=\left\{z_S\right\}$ on $V_1 \times V_2 \times V_3$ by letting $z_S=\left((4r+2)(2r+1)\right)^{-1}$ if $S=(a, b, c)$ with $a+b+c$ even and $z_S=0$ otherwise. Then $Z$ is a 3-stochastic weight on $V_1 \times V_2 \times V_3$ while 
$P(Z)=0$, since the sum of all integers in $V_1$, $V_2$ and $V_3$ is odd. Hence the permanent of a 3-stochastic 
$m \times m \times m$ tensor can be zero even when the maximum entry of the tensor is $O\left(m^{-2}\right)$.
This example was constructed in a conversation with Jeff Kahn.

Summarizing, for general $k$-stochastic weights $Z$ there is no a priori non-zero lower bound for the partition function. If, however, we require $Z$ to be $\alpha$-balanced for any fixed 
$\alpha \geq 1$, the lower bound jumps to within a polynomial in $m$ factor of the upper bound.

We note that there are extensions of the Bregman-Minc bound to hafnians \cite{AF08}
and to higher-dimensional permanents \cite{D87a}. Lower bounds appear to be harder to come by, see \cite{E+10} for the recent proof of the Lov\'asz-Plummer conjecture, which states that the number of perfect 
matchings in a bridgeless 3-regular graph is exponentially large in the number of vertices of the graph, and 
\cite{F11b} and \cite{Ba11} for related developments. 
\bigskip
If $H$ is a complete $k$-uniform hypergraph or a complete $k$-partite hypergraph, one can estimate $P_H(W)$ for any balanced but not necessarily $k$-stochastic weight $W$ using scaling.

\subhead (1.5) Scaling \endsubhead Let $W=\left\{w_S\right\}$ be a weight on the edges of a $k$-uniform hypergraph 
$H$
with a vertex set $V$, where $|V|=km$. Let $\left\{\lambda_v>0: \ v \in V\right\}$ be reals. We say that 
a weight $Z=\left\{z_S\right\}$ on the hypergraph $H$ is obtained from $W$ by {\it scaling} if 
$$z_S= \left( \prod_{v \in S} \lambda_v \right) w_S \quad \text{for all} \quad S \in H.$$
It is easy to see that 
$$P_H(Z) = \left( \prod_{v \in V} \lambda_v \right) P_H(W).$$
It turns out that any positive weight $W$ on a complete $k$-uniform hypergraph or a complete 
$k$-partite hypergraph can be scaled to a $k$-stochastic weight $Z$
(cf., for example, \cite{F11a} and Section 3 below). We show that the $k$-stochastic scaling of an $\alpha$-balanced 
weight is $\alpha^{k+1}$-balanced and obtain the following result.

\proclaim{(1.6) Theorem} Let us fix an integer $k >1$ and a real $\alpha \geq 1$. Then there exists a
real $\gamma=\gamma(k, \alpha)>0$ such that the following holds.

Let $H$ be a complete $k$-uniform hypergraph  or a complete $k$-partite hypergraph with 
$km$ vertices and let $W=\left\{w_S: \ S \in H \right\}$ be an $\alpha$-balanced weight on $H$.
Let us consider the function  
$$f_W(X) =\sum_{S \in H} x_S \ln {x_S \over w_S}$$
for a weight $X$ on $H$. Let $\Omega_k(H)$ be the set of all $k$-stochastic weights 
on $H$ and let 
$$\zeta=\min_{X \in \Omega_k(H)} f_W(X).$$
Then
$$e^{-\zeta -m(k-1)} m^{-\gamma} \ \leq \ P_H(W) \ \leq \ e^{-\zeta-m(k-1)} m^{\gamma}.$$
\endproclaim

More precisely, we prove that under conditions of Theorem 1.3 and assuming, additionally, that $\alpha^{k+1}>2$, we have 
$$\epsilon_1 m^{-\gamma_1} e^{-\zeta-m(k-1)}  \ \leq \ P_H(W) \ \leq \ \epsilon_2 m^{\gamma_2} e^{-\zeta-m(k-1)},$$
where $\gamma_1, \gamma_2, \epsilon_1, \epsilon_2$ are defined by (1.3.1).

The set $\Omega_k(H)$ is naturally identified with a convex polytope in ${\Bbb R}^{H}$. Function $f$ is strictly 
convex and hence the optimization problem of computing $\zeta$ can be solved efficiently (in polynomial time)
by interior point methods, see \cite{NN94}. Thus Theorem 1.6 allows us to estimate the partition function of 
an $\alpha$-balanced weight (for any $\alpha \geq 1$, fixed in advance) within a polynomial in $m$ factor. 

\subhead (1.7) A probabilistic interpretation \endsubhead Let us fix $k>1$ and 
let $H$ be either a complete $k$-uniform hypergraph or a complete $k$-partite hypergraph 
with a set $V$ of $|V|=km$ vertices.
Let us fix $\alpha \geq 1$ and let $W=\left\{w_S: \quad S \in H\right\}$ be an 
$\alpha$-balanced weight on $H$. Let $|H|$ denote the number of edges of 
hypergraph $H$.
 Let us assume that 
$$\sum_{S \in H} w_S =m,$$
in which case 
$${m \over \alpha |H|} \ \leq \ w_S \ \leq \ {\alpha m \over |H|} \quad 
\text{for all} \quad S \in H.$$
In particular, for all sufficiently large $m$ we have $w_S < 1$ for all $S \in H$, so  
we can introduce independent random Bernoulli 
variables $X_S$ indexed by the edges $S \in H$, where 
$$\PP\left(X_S=1\right) =w_S \quad \text{and} \quad  \PP\left(X_S=0\right)=1-w_S.$$
For each vertex $v \in V$ let us define a random variable 
$$Y_v =\sum \Sb S \in H \\ S \ni v  \endSb X_S.$$
It is not hard show that 
$$P_H(W) =\exp\left\{m + O\left({1 \over m^{k-2} }\right) \right\} 
\PP\Bigl(Y_v = 1 \quad \text{for all} \quad v \in V \Bigr).$$
For large $m$, the distribution of each random variable $Y_v$ is approximately Poisson with 
$$\EE Y_v=\mu_v \quad \text{where} \quad \mu_v= \sum\Sb S \in H  \\ S \ni v \endSb w_S,$$
so 
$$\PP\left(Y_v=1 \right) \approx \mu_v e^{-\mu_v}.$$
The probability of $Y_v=1$ is maximized when $\mu_v=1$, and when $W$ is $k$-stochastic, 
the probabilities of $Y_v=1$ are maximized simultaneously for all $v \in V$, so that
$$\PP\left(Y_v=1\right) \approx e^{-1} \quad \text{for all} \quad v \in V.$$
Theorem 1.3 implies that in this case the events $Y_v=1$ behave as if they were (almost) independent, so that 
$$\PP\Bigl(Y_v = 1 \quad \text{for all} \quad v \in V \Bigr) \approx  e^{-km}$$
up to a polynomial in $m$ factor.
\bigskip
In Section 2, we discuss some combinatorial and algorithmic applications of Theorems 1.3 and 1.6. Namely, we 
present a simple polynomial time algorithm to distinguish hypergraphs having sufficiently many perfect 
matchings from hypergraphs that do not have nearly perfect matchings. We also prove a lower bound for the 
number of nearly perfect matchings in regular hypergraphs.

In the rest of the paper we prove Theorems 1.3 and 1.6.

In Section 3, we review some results about scaling. The results are not new, but we nevertheless provide 
proofs for completeness. In Section 4, we prove two crucial lemmas about scaling of $\alpha$-balanced weights.
In Section 5 we complete the proofs of Theorems 1.3 and 1.6.

Scaling was used in \cite{L+00} to efficiently estimate permanents of non-negative matrices.   

\subhead (1.8) Notation \endsubhead As usual, for two functions $f$ and $g$, where $g$ is 
non-negative, we say that $f=O(g)$ if $|f| \leq \gamma g$ for some constant $\gamma >0$.
We will allow our constants $\gamma$ to depend only on the dimension $k$ of the 
hypergraph and the parameter $\alpha \geq 1$ in the definition of an $\alpha$-balanced 
weight in Section 1.2.

\head 2. Combinatorial applications \endhead

 Let us fix an integer $k>1$
and let $H$ be a $k$-uniform hypergraph with $km$ vertices.
As is known \cite{Va79}, the problem of counting perfect matchings 
in $H$ is $\#$P-hard. For $k=2$ there is a classical polynomial time algorithm to check whether $H$ has a perfect matching (see \cite{LP09}) and a fully polynomial randomized approximation scheme is known for counting perfect matchings if $H$ is 
bipartite \cite{J+04}. For $k>2$ finding if there is a perfect matching in $H$ is an NP-complete problem 
even when $H$ is $k$-partite \cite{Ka72}.

Theorem 1.6 allows us to distinguish in polynomial time between hypergraphs that have sufficiently many perfect matchings 
and hypergraphs that do not have nearly perfect matchings.

In this section, we let 
$$\Phi_k(m)={(km)! \over (k!)^m m!}$$
be the number of perfect matchings in a complete $k$-uniform hypergraph with $km$ vertices.
\subhead (2.1) Testing hypergraphs \endsubhead 
Let us fix integer $k>1$ and positive real $\delta \leq 1$ and $\beta < 1$. 

We consider the following algorithm.
\bigskip
{\bf Input:} A $k$-uniform hypergraph $H$,  defined by the list of its edges, 
with a set $V$ of $km$ vertices.
\medskip
{\bf Output:} At least  one of the following two (not mutually exclusive) conclusions:
\medskip
(a) The hypergraph $H$ contains a matching with at least $\beta m$ edges.
\smallskip
(b) The hypergraph $H$ contains at most $\delta^m \Phi_k(m)$ perfect matchings.
\medskip
{\bf Algorithm:} Let  
$$\epsilon ={1 \over 2} \delta^{1/(1-\beta)}.$$
 Let us define a weight 
$W=\left\{w_S\right\}$ on the complete $k$-uniform hypergraph ${V \choose k}$ as follows:
$$w_S=\cases 1 & \text{if \ } S \in H \\ \epsilon &\text{if \ } S \notin H. \endcases \tag2.1.1$$
The weight $W$ is $\epsilon^{-1}$-balanced and we apply the algorithm of 
Theorem 1.6 to compute in polynomial in $m$ time a number $\eta$ such that 
$$\eta \cdot m^{-\gamma} \ \leq \ P(W) \ \leq \ \eta \cdot m^{\gamma}$$
for some $\gamma=\gamma(\delta, \beta)>0$. 

If $m=1$ or if
$${m \over \ln m} \ \leq \ {2 \gamma \over (1-\beta) \ln 2} \tag2.1.2$$
we check by direct enumeration whether (a) or (b) hold. Since $k$, $\beta$ and $\delta$ are fixed 
in advance, this requires only a constant time.

If  (2.1.2) does not hold, we output conclusion (a) if
 $\eta \cdot m^{\gamma} > \delta^m \Phi_k(m)$ and conclusion (b) if 
 $\eta \cdot m^{\gamma} \leq \delta^m \Phi_k(m)$.
\bigskip
It is not hard to see that the algorithm is indeed correct.  
If $\eta \cdot m^{\gamma} \leq \delta^m \Phi_k(m)$ then $P(W) \leq \delta^m \Phi_k(m)$ and 
$H$ necessarily contains not more than $\delta^m \Phi_k(m)$ perfect matchings.
If $\eta \cdot m^{\gamma} > \delta^m \Phi_k(m)$ then, assuming that (2.1.2) does not hold,
we conclude that 
$$P(W) \ \geq \ {\delta^m \over m^{2 \gamma}} \Phi_k(m) \ > \  {\delta^m \over 2^{(1-\beta)m}} \Phi_k(m)  
\ = \ \epsilon^{(1-\beta) m} \Phi_k(m),$$
from which it follows that $H$ contains a matching with not fewer than $\beta m$ edges.

In particular, confronted with two hypergraphs on $km$ vertices, one of which contains more than $\delta^m \Phi_k(m)$ 
perfect matchings and the other with no matchings of size $\beta m$ or bigger, the algorithm will be able to decide 
which is which. It will necessarily output a) in the former case and b) in the latter.

\definition{(2.2) Definition} 
A $k$-uniform hypergraph $H$ is called {\it $d$-regular} if every vertex of $H$ is contained in exactly $d$ 
edges of $H$.
\enddefinition
For example, a complete $k$-uniform hypergraph with $km$ vertices is $d$-regular for 
$d={km-1 \choose k-1}$. The existence
of a perfect or nearly perfect matching in $d$-regular hypergraphs was extensively studied, see, for example, \cite{Vu00}
and references therein.  As a corollary of Theorem 1.3, we obtain the following estimate for the number of 
nearly perfect matchings in a regular hypergraph (see also \cite{C+91} for some related estimates).

\proclaim{(2.3) Theorem} Let us fix an integer $k>1$ and $0< \alpha, \beta <1$.
Then there exists a positive integer $m_0=m_0(k, \alpha, \beta)$ such that the following holds.

Let $H$ be a $k$-uniform $d$-regular hypergraph with $km$ vertices where 
$$d \ \geq \ \alpha {km-1 \choose k-1} \quad \text{and} \quad m \geq m_0.$$ 
Then for every positive integer $s \leq \beta m$ the hypergraph $H$
contains at least  
$$ \alpha^m {\Phi_k(m) \over \Phi_k(m-s)}$$
matchings of size $s$.
\endproclaim
\demo{Proof} All implied constants in the ``$O$'' notation below may depend on $k, \alpha$ and $\beta$ only.

 Let $V$ be the set of vertices of a $k$-uniform $d$-regular hypergraph, $|V|=km$.
Let us choose $0 < \epsilon < 1$ such that 
$$\epsilon^{1-\beta} \ < \ \alpha + \epsilon(1-\alpha) \tag2.3.1$$ 
and let us define a weight $W=\left\{w_S\right\}$ on the complete $k$-uniform hypergraph ${V \choose k}$ by (2.1.1).
Then 
$$\sum \Sb S \in {V \choose k} \\ S \ni v \endSb w_S =(1-\epsilon) d + \epsilon {km-1 \choose k-1} \quad 
\text{for all} \quad v \in V.$$
It follows from Theorem 1.3 and scaling that
$$\aligned P(W)\ \geq \ &\left((1-\epsilon) d +\epsilon {km-1 \choose k-1}\right)^m e^{-m(k-1)} {1 \over m^{O(1)}} \\
 \geq \ &\left(\alpha+\epsilon(1-\alpha) \right)^m {km -1 \choose k-1}^m e^{-m(k-1)} {1 \over m^{O(1)}}. \endaligned
 \tag2.3.2$$
We note that 
$$\aligned {{km -1\choose k-1}^m \over \Phi_k(m)} =
&{\bigl((km-1)!\bigr)^m (k!)^m m! \over \bigl((k-1)!\bigr)^m \bigl((km-k)!\bigr)^m(km)!}=
{\bigl((km-1)!\bigr)^m k^m m! \over  \bigl((km-k)!\bigr)^m(km)!}\\=
&{\bigl((km)!\bigr)^m k^m m! \over  (km)^m \bigl((km-k)!\bigr)^m(km)!}= 
{\bigl((km)!\bigr)^m  m! \over  \bigl((km-k)!\bigr)^m (km)! m^m}\\
=&{\left((km)(km-1) \cdots (km-k+1) \over (km)^k \right)^m} \cdot   {(km)^{km} m! \over (km)! m^m}. \endaligned
\tag2.3.3$$
Since $\ln (1-x) \geq -2x$ for $0 \leq x \leq 0.5$, we conclude that 
$${\left((km)(km-1) \cdots (km-k+1) \over (km)^k \right)^m} =\exp\left\{m \sum_{i=1}^{k-1} \ln \left(1 -{i \over km}\right)\right\}
\ \geq \ e^{-k+1}.$$
Using Stirling's formula, we conclude from (2.3.3) and (2.3.2) that 
$$P(W) \ \geq \ \left(\alpha+\epsilon(1-\alpha) \right)^m \Phi_k(m) {1 \over m^{O(1)}}. \tag2.3.4$$

If a perfect matching in ${V \choose k}$ contains 
fewer than $s$
edges of $H$ then the contribution of the corresponding term to $P(W)$ is less than $\epsilon^{m-s}$. Since every matching in 
$H$ of size $s$ can be appended to a perfect matching in ${V \choose k}$ in $\Phi_k(m-s)$ ways, we conclude 
that the number of matchings in $H$ of size $s$ is at least 
$${P(W) - \epsilon^{m-s} \Phi_k(m) \over  \Phi_k(m-s)} \ \geq \ {P(W) - \epsilon^{(1-\beta)m} \Phi_k(m) \over  \Phi_k(m-s)} .$$
The proof now follows from (2.3.4) and (2.3.1).
{\hfill \hfill \hfill} \qed
\enddemo

\head 3. General results on scaling \endhead

In this section, we summarize some results on scaling which we need for the proofs of Theorems 1.3 and 1.6.

\proclaim{(3.1) Theorem} Let $H$ be a $k$-uniform hypergraph with a set $V$ of $|V|=km$ 
vertices and let $\Omega_k(H)$ be the set of all $k$-stochastic weights on $H$. Suppose 
that the set $\Omega_k(H)$ has a non-empty relative interior, that is contains a positive weight $Y \in \Omega_k(H)$.

For a positive weight $W=\left\{w_S: \ S \in H\right\}$ on $H$, let us define a 
function $f_W: \Omega_k(H) \longrightarrow {\Bbb R}$ by
$$f_W(X)=\sum_{S \in H} x_S \ln {x_S \over w_S} \quad \text{for} \quad X \in \Omega_k(H),
\ X=\left\{x_S: \ S \in H \right\}.$$
Then function $f_W$ attains its minimum on $\Omega_k(H)$ at a unique weight \break
$Z=\left\{z_S: \ S \in H \right\}$. We have $z_S >0$ for all $S \in H$ and 
there exist real \break $\lambda_v>0: \ v \in V$ such that 
$$z_S =\left(\prod_{v \in S} \lambda_v\right) w_S \quad \text{for all} \quad S \in H. \tag3.1.1$$ 
We have
$$f_W(Z) =\sum_{v \in V}  \ln  \lambda_v.$$
Furthermore, if $\lambda_v>0: v \in V$ are reals such that weight $Z$ 
defined by (3.1.1) is $k$-stochastic, then $Z$ is the minimum point of $f_W$ on $\Omega_k(H)$.
\endproclaim
\demo{Proof} First, we observe that function $f_W$ is strictly convex, so its minimum on the convex set 
$\Omega_k(H)$ is unique. Next,
$${\partial \over \partial x_S} f_W(X) =\ln {x_S \over w_S} +1, \tag3.1.2$$
which is finite if $x_S >0$ and is $-\infty$ if $x_S=0$ (we consider the right derivative in this case).
If $z_S=0$ for some $S$ then for a sufficiently small $\epsilon >0$ we have 
$$f_W\bigl((1-\epsilon)Z + \epsilon Y \bigr)  < f_W(Z),$$
which is a contradiction. Hence $z_S>0$ for all $S \in H$. 

Since the minimum point $Z$ lies in the relative interior of $\Omega_k(H)$, considered 
as a convex polyhedron in ${\Bbb R}^H$, the Lagrange multiplier condition 
implies that 
there exist real $\mu_v:\ v \in V$ such that 
$$\ln {z_S \over w_S} =\sum_{v \in S} \mu_v \quad \text{for all} \quad S \in H. \tag3.1.3$$
Hence, letting $\lambda_v = e^{\mu_v}$ for $v \in V$, we obtain
$$z_S =\left( \prod_{v \in S} \lambda_v \right) w_S \quad \text{for all} \quad S \in H.$$
Now, 
$$f_W(Z)=\sum_{S \in H} z_S \left( \sum_{v \in S} \ln \lambda_v \right) =
\sum_{v \in V} \ln \lambda_v \left( \sum \Sb S \in H \\ S \ni v \endSb z_S \right) =
\sum_{v \in V} \ln \lambda_v,$$
as desired.

If (3.1.1) holds for some $\lambda_v>0$ and $k$-stochastic $Z=\left\{z_S\right\}$, then 
(3.1.3) holds with $\mu_v =\ln \lambda_v$ and by (3.1.2) we conclude that $Z$ is a critical 
point of $f_W$ in the relative interior of $\Omega_k(H)$. Since $f_W$ is strictly convex, $Z$ must 
be the minimum point of $f_W$ on $\Omega_k(H)$.
{\hfill \hfill \hfill} \qed
\enddemo
Theorem 3.1 implies that any positive weight $W$ on a hypergraph $H$ having a positive 
$k$-stochastic weight can be scaled uniquely to a $k$-stochastic weight $Z$, in which case we have
$P_H(W)=\exp\{ - f_W(Z) \} P_H(Z)$. Scaling factors $\left\{ \lambda_v>0: \ v \in V \right\}$, however, do not have to be 
unique, as the example of a complete $k$-partite hypergraph readily shows (although in the case of the 
complete $k$-uniform hypergraph the scaling factors are unique).
We note that if $H$ is the complete $k$-uniform 
hypergraph or the complete $k$-partite hypergraph then there is a positive 
$k$-stochastic weight $Y=\left\{ y_S: \ S \in H \right\}$ on $H$.
In the former case 
we can choose 
$$y_S={km-1 \choose k-1}^{-1} \quad \text{for all} \quad S \in H,$$
while in the latter case we can choose
$$y_S=m^{-k+1} \quad \text{for all} \quad S \in H.$$
We need a dual description of the scaling factors $\lambda_v$.
\proclaim{(3.2) Theorem} Let $H$ be a $k$-uniform hypergraph with a set $V$ of $|V|=km$ 
vertices and let $W=\left\{w_S: \ S \in H \right\}$ be a positive weight on $H$. 
Let $\lambda_v >0: \ v \in V$ be reals such that the weight $Z=\left\{z_S\right\}$ defined by 
$$z_S =\left(\prod_{v \in S} \lambda_v \right) w_S \quad \text{for all} \quad S \in H$$
is $k$-stochastic.

Let us define a set $C(W) \subset {\Bbb R}^V$ by 
$$C(W)=\left\{ \left(x_v, \ v \in V \right): \quad \sum_{S \in H}
w_S \exp\left\{ \sum_{v \in S} x_v \right\} \leq m \right\}.$$
Then the point $\left(\mu_v: \ v \in V\right)$, where $\mu_v=\ln \lambda_v$ for all $v \in V$, is 
a maximum point of the linear function $\sum_{v \in V} x_v$ on 
$C(W)$.
\endproclaim
\demo{Proof} Since weight $Z$ is $k$-stochastic, we have 
$$\sum \Sb S \in H \\  S \ni u \endSb w_S \exp\left\{\sum_{v \in S} \mu_v \right\} 
=1 \quad \text{for all} \quad u \in V,$$ 
which means that $\left(\mu_v: \ v \in V \right)$ is a critical point of the linear function 
$\sum_{v \in V} x_v$ on the smooth surface defined in ${\Bbb R}^V$ by the equation 
$$\sum_{S \in H} w_S \exp\left\{ \sum_{v \in S} x_v \right\} =m.$$
The set $C(W)$ is convex and hence $\left(\mu_v: \ v \in V \right)$ has to be an extremum point 
of function $\sum_{v \in V} x_v$ on $C(W)$. Moreover, it has to be a maximum point 
since the function is unbounded from below on $C(W)$. 
{\hfill \hfill \hfill} \qed
\enddemo

\head 4. Scaling balanced weights \endhead 

Our proof of Theorem 1.3 is based on two lemmas. 

\proclaim{(4.1) Lemma} Let us fix an integer $k>1$ and a real $\alpha  \geq 1$
and let $H$ be a complete $k$-uniform hypergraph or a complete $k$-partite hypergraph.
If $W=\left\{w_S\right\}$ is an $\alpha$-balanced weight on 
$H$ and if $Z=\left\{z_S\right\}$ is the $k$-stochastic weight obtained 
from $W$ by scaling, then $Z$ is $\alpha^{k+1}$-balanced.
\endproclaim 
\demo{Proof} Let $V$ be the set of vertices of hypergraph $H$. Without loss of generality, we assume
that $|V|>k$. For a subset $X \subset V$, we denote 
by 
$$H_X=\left\{S \in H: \quad  S  \supset X\right\}$$
the set of edges of $H$ containing $X$. 
Let $\left\{\lambda_v>0: \ v \in V\right\}$ be scaling factors so that 
$$z_S=\left( \prod_{v \in S} \lambda_v \right) w_S \quad \text{for all} \quad S \in H.\tag4.1.1$$

Suppose first that $H=V_1 \times \ldots \times V_k$ is a complete $k$-partite hypergraph, so
$V=V_1 \cup \ldots \cup V_k$ and $\left|V_1\right| =\ldots =\left| V_k \right|$.
 For every $i=1, \ldots, k$ and for every pair of vertices $v, u \in V_i$ we have 
$$\sum_{S \in H_{\{v\}}} z_S = \sum_{S \in H_{\{u\}}} z_S =1. \tag4.1.2$$
Let us consider the bijection $\phi: H_{\{v\}}  \longrightarrow H_{\{u\}}$ defined by 
$$ \phi(S)= S \cup\{u\} \setminus \{v\}. \tag4.1.3$$
By (4.1.1) we have
$${z_{\phi(S)}\over z_{S}} ={\lambda_u \over \lambda_v} \cdot {w_{\phi(S)} \over w_S}. \tag4.1.4$$
Since weight $W$ is $\alpha$-balanced, in view of (4.1.2) we conclude that 
$${1 \over \alpha} \ \leq \ {\lambda_u \over \lambda_v} \ \leq \ \alpha, \tag4.1.5$$
which proves that $Z$ is $\alpha^{k+1}$-balanced. 

Suppose now that $H={V \choose k}$ is a complete $k$-uniform hypergraph.
Then for any two distinct vertices $u, v \in V$ we have 
$$\sum_{S \in H_v \setminus H_{\{u, v\}}} z_S = \sum_{S \in H_u \setminus H_{\{u, v\}}} z_S =
1-\sum_{S \in H_{\{u, v\}}} z_S >0. \tag4.1.6$$
Let us consider the bijection $\phi: H_v \setminus H_{\{u, v\}} \longrightarrow H_u \setminus H_{\{u, v\}}$
defined by (4.1.3). From (4.1.1) we deduce that (4.1.4) holds and in view of (4.1.6) we conclude that (4.1.5) 
follows. Since weight $W$ is $\alpha$-balanced, (4.1.5) implies that $Z$ is $\alpha^{k+1}$-balanced. 
{\hfill \hfill \hfill} \qed
\enddemo

The second lemma asserts that if we scale to a $k$-stochastic weight a weight which is 
already sufficiently close to being $k$-stochastic, then the product of the scaling factors is close to 1.

\proclaim{(4.2) Lemma} Let us fix an integer $k>1$ and real $\alpha \geq 1$ and $\delta >0$. 
Then there exist integer $m_0=m_0(k, \alpha, \delta) >0$ and real $\beta=\beta(k, \alpha, \delta) > 0$
such that the following holds.

Suppose that $H$ is a complete $k$-uniform hypergraph or a complete $k$-partite hypergraph with a set 
$V$ of vertices, where $|V|=km$ and $m \geq m_0$.
Suppose that $W=\left\{w_S \right\}$ is an $\alpha$-balanced weight on $H$, that
$$\sum_{S \in H} w_S =m$$ and that
$$\left| 1- \sum \Sb S \in H \\ S \ni v 
\endSb  w_S  \right|  \ \leq \  {\delta \over m} 
\quad \text{for all} \quad v \in V.$$
Let $\lambda_v>0: \ v \in V $ be reals such that 
weight $Z=\left\{ z_S\right\}$ defined by 
$$z_S =\left( \prod_{v \in S} \lambda_v \right) w_S \quad \text{for all} \quad S \in H$$
is $k$-stochastic.
Then 
$$0\  \leq \ \sum_{v \in V} \ln \lambda_v \ \leq \ {\beta \over m}.$$
One can choose $\beta=\alpha \delta^2 (k+1)^2$ and  
$m_0  =  \max\bigl\{1+\lceil \alpha \delta k\rceil, \ k \bigr\}$.
\endproclaim
\demo{Proof} 
We note that the point $\bigl(x_v=0: \ v \in V \bigr)$ belongs to the set $C(W)$ of Theorem 3.2, and 
so by Theorem 3.2 we have
$$\sum_{v \in V}  \ln \lambda_v \ \geq \ \sum_{v \in V} x_v = 0.$$
Let us define
$$\delta_v = 1-\sum \Sb S \in H \\ S \ni v \endSb w_S \quad \text{for} \quad v \in V.$$
Then
$$\sum_{v \in V} \delta_v = \sum_{v \in V} \left(1- \sum \Sb S \in H \\ S \ni v 
\endSb w_S  \right)
=km- k \sum_{S \in H} w_S  =0. \tag4.2.1$$
In addition, if $H=V_1 \times \ldots \times V_k$ is the complete $k$-partite graph, where 
$V=V_1 \cup \ldots \cup V_k$ and $\left| V_1\right| =\ldots =\left|V_k\right|=m$, for every 
$i=1, \ldots, k$ we have 
$$\sum_{v \in V_i} \delta_v =\sum_{v \in V_i} \left(1- \sum \Sb S \in H \\ S \ni v 
\endSb w_S  \right)=m- \sum_{S \in H} w_S =0. \tag4.2.2$$
We define numbers $\left\{ \rho_S: \ S \in H \right\}$ as follows.
If $H={V \choose k}$ is the complete $k$-uniform hypergraph, 
we define 
$$\rho_S ={km-2 \choose k-1}^{-1} \sum_{v \in S} \delta_v \quad \text{for all} \quad 
S \in H.$$
If $H=V_1 \times \ldots \times V_k$ is a complete $k$-partite graph, we define 
$$\rho_S = {1 \over m^{k-1}} \sum_{v \in S} \delta_v \quad \text{for all} \quad S \in H.$$
We claim that 
$$\sum \Sb S \in H \\ S \ni v \endSb \rho_S =\delta_v \quad \text{for all} \quad v \in V. \tag4.2.3$$
Indeed, if $H={V \choose k}$ then using (4.2.1) we obtain
$$\split  \sum \Sb S \in H \\ S \ni v \endSb \rho_S  = 
&{{km-1 \choose k-1} \over {km-2 \choose k-1}}  \delta_v + 
{{km-2 \choose k-2} \over{km-2 \choose k-1}} \sum_{u \in V \setminus \{v\}} \delta_u \\
=&{ {km-1 \choose k-1} -{km-2 \choose k-2} \over {km-2 \choose k-1}} \delta_v  \\
  =& \delta_v 
\endsplit$$
and if $H=V_1 \times \ldots \times V_k$ then using (4.2.2)
we obtain that for all $i=1, \ldots, k$ and for all $v \in V_i$ we have 
$$\sum \Sb S \in H \\ S \ni v \endSb  \rho_S =\delta_v +{m^{k-2}\over m^{k-1}} \sum_{u \in V\setminus V_i} \delta_u=
\delta_v. $$
In either case, (4.2.3) holds. In addition, from (4.2.1)
$$\sum_{S \in H} \rho_S = {1 \over k} \sum_{v \in V} \sum \Sb S \in H \\ S \ni v \endSb \rho_S =
{1 \over k} \sum_{v \in V} \delta_v =0. \tag4.2.4$$
Let us define
$$x_S = w_S + \rho_S \quad \text{for all} \quad S \in H.$$
Then, from (4.2.3) we have 
$$\sum \Sb S \in H \\ S \ni v \endSb x_S =1 \quad \text{for all} \quad v \in V. \tag4.2.5$$
Since weight $W$ is $\alpha$-balanced, for all $S \in H$ we have 
$$w_S \ \geq \ {km \choose k}^{-1}{m \over \alpha} \quad \text{when} \quad H={V \choose k}
\tag4.2.6$$
and
$$w_S \ \geq \ {1 \over \alpha m^{k-1}} \quad \text{when} \quad 
H=V_1 \times \ldots \times V_k.
\tag4.2.7$$
On the other hand, for all $S \in H$ we have 
$$\left| \rho_S \right| \ \leq \  {km-2 \choose k-1}^{-1} {\delta k \over m } \quad \text{when} 
\quad H= {V \choose k} \tag4.2.8$$
and 
$$\left| \rho_S \right| \ \leq \ {\delta k \over m^k} \quad \text{when} \quad H=V_1 \times \ldots \times V_k. \tag4.2.9$$
From (4.2.6) and (4.2.8)
we conclude that if $H={V \choose k}$ then
$$x_S \ \geq \ 0 \quad \text{for all} \quad S \in H \quad \text{provided} 
\quad {m (m-1) \over km-1} \ \geq \ \alpha \delta $$
whereas from (4.2.7) and (4.2.9) we conclude that if $H=V_1 \times \ldots \times V_k$ then 
$$x_S \ \geq \ 0 \quad \text{for all} \quad S \in H \quad \text{provided} 
\quad m \ \geq \ \alpha \delta k. $$
In either case, $X=\left\{x_S \right\}$ is a $k$-stochastic weight on $H$ provided 
$m \geq \alpha \delta k +1$.
Using (4.2.4),  we conclude from Theorem 3.1 that for $m \geq \alpha \delta k +1$
we have
$$\aligned \sum_{v \in V} \ln \lambda_v \ \leq \ &\sum_{S \in H} x_S \ln {x_S \over w_S} = 
\sum_{S \in H} \left(w_S + \rho_S\right) \ln {w_S + \rho_S \over w_S} \\ = & \sum_{S \in H}
\left(w_S +\rho_S \right) \ln \left(1 + {\rho_S \over w_S}\right) 
\ \leq \ \sum_{S \in H} \left(w_S +\rho_S\right) {\rho_S \over w_S} \\= 
&\sum_{S \in H} {\rho_S^2 \over w_S}.
\endaligned \tag4.2.10$$
From (4.2.6) and (4.2.8), we conclude that in the case of $H={V \choose k}$ sum (4.2.10) does not exceed
$${\alpha \delta^2 k^2{km \choose k}^2  \over m^3 {km-2 \choose k-1}^2 }={\alpha \delta^2 (km-1)^2 \over m (m-1)^2}$$
whereas from (4.2.7) and (4.2.9) we conclude that in the case of $H=V_1 \times \ldots \times V_k$ 
sum (4.2.10) does not exceed 
$${\alpha \delta^2 k^2 m^{k-1} \over m^{2k}} m^k = {\alpha \delta^2 k^2 \over m}.$$
In either case, sum (4.2.10) does not exceed $\alpha \delta^2(k+1)^2/m$ as long as $m \geq k$.
{\hfill \hfill \hfill} \qed
\enddemo

\head 5. Proofs of Theorems 1.3 and 1.6 \endhead

Our approach is somewhat similar to Bregman's original approach \cite{Br73} combining scaling and induction 
to obtain upper bounds on permanents. Before giving a formal proof, we illustrate the idea of the proof 
by sketching it in the more familiar case of permanents, that is when $k=2$ and 
the underlying graph is bipartite.

All implied constants in the ``$O$'' notation in this section may depend on $k$ and $\alpha$ only. 

\subhead (5.1) The idea of the proof \endsubhead 
Let us fix $\alpha \geq 1$. Let $A=\left(a_{ij}\right)$ be an $\alpha$-balanced $m \times m$ doubly stochastic matrix.
Our goal is to prove that $\per A=e^{-m} m^{O(1)}$, or, equivalently, that 
$$\per A = \exp\left\{-m+O\left( \sum_{j=1}^m {1 \over j}\right)\right\}. \tag5.1.1$$
We proceed by induction on $m$.

Using the first row expansion, we can write 
$$\per A =\sum_{j=1}^m a_{1j} \per \widehat{A}_{j}, \tag5.1.2$$
where $\widehat{A}_j$ is the $(m-1)\times (m-1)$ matrix obtained from $A$ by crossing out the first row and $j$-th 
column. We have $a_{1j} \leq \alpha/m$ for all $j$ and, using that $A$ is doubly 
stochastic, we conclude that the sum of 
$\sigma_j$ of all entries of $\widehat{A}_j$ satisfies 
$$ m-2 \ \leq \ \sigma_j \ \leq \ m-2 + {\alpha \over m}. \tag5.1.3$$
Let us define $(m-1)\times (m-1)$ matrices $B_j$ by 
$$B_j = {m-1 \over \sigma_j} \widehat{A}_j.$$
Hence the sum of all entries of $B_j$ is $m-1$ and by (5.1.3) we have 
$$\per \widehat{A}_j = \left({\sigma_j \over m-1}\right)^{m-1} \per B_j = \exp\left\{ -1+ O\left({1 \over m}\right)\right\}
\per B_j. \tag5.1.4$$

Matrices $B_j$ are not necessarily doubly stochastic, but they are reasonably close to doubly stochastic, since all 
the row and column sums of $B_j$ are $1+O\left(m^{-1}\right)$. Let $C_j$ be the $(m-1) \times (m-1)$ doubly stochastic 
matrix obtained from $B_j$ by scaling. By Lemma 4.2 we have 
$$\per B_j = \exp\left\{ O\left({1 \over m}\right)\right\} \per C_j. \tag5.1.5$$
Combining (5.1.2) and (5.1.3)--(5.1.5), we conclude that 
$$\per A = \exp\left\{ -1 + O\left({1 \over m}\right)\right\} \sum_{j=1}^m a_{1j} \per C_j, \tag5.1.6$$
where 
$$\sum_{j=1}^m a_{1j} =1 \quad \text{and} \quad a_{1j} \geq 0 \quad \text{for all} \quad j=1, \ldots, m.$$
Up until this point, we did not really use the condition that $A$ is $\alpha$-balanced, we used only that 
the entries of $A$ are uniformly small, of the order of $O\left(m^{-1}\right)$. To proceed with the induction, we have 
to show that the entries of the doubly stochastic matrices $C_j$ in (5.1.6) and all other doubly stochastic matrices 
obtained by iterating the recursion are also uniformly small. Now we observe that $C_j$ is obtained by scaling of 
$\widehat{A}_j$ and hence by Lemma 4.1 is $\alpha^3$-balanced. Similarly, as we iterate recursion (5.1.6), 
the doubly stochastic matrices that we obtain are $\alpha^3$-balanced, since they are obtained by scaling 
from some submatrices of an $\alpha$-balanced matrix $A$. This allows us to use (5.1.6) in the induction step 
to obtain (5.1.1).

Permanents of $\alpha$-balanced matrices are studied in \cite{F+04} and \cite{CV09}.

\subhead (5.2) Proof of Theorem 1.3 \endsubhead Without loss of generality, we assume that $\alpha^{k+1} > 2$.

Let $H$ be either a complete $k$-uniform hypergraph ${V \choose k}$ or a complete $k$-partite hypergraph 
$V_1 \times \ldots \times V_k$ with a set 
$V$ of $|V|=km$ vertices.
Let $Z=\left\{z_S \right\}$ be a $k$-stochastic $\alpha$-balanced weight on $H$. 

If $H={V \choose k}$ and
$U \subset V$ is a subset such that $|U|=kl$ for some integer $l \geq1$, we consider the induced hypergraph 
$H|U$ consisting of the edges $S \in H$ such that $S \subset U$. Hence $H|U={U \choose k}$ is the complete $k$-uniform hypergraph with the set $U$ of vertices. 

Similarly, if $H=V_1 \times \ldots \times V_k$ and if 
$\left|U \cap V_1 \right| =\ldots =\left|U \cap V_k \right| =l$ for some integer $l \geq1$, we consider the 
restriction $H|U$ consisting of the edges $S \in H$ such that $S \subset U$. In this case, $H|U$ is a uniform $k$-partite 
graph with the set $U$ of vertices, $H|U=U_1 \times \ldots \times U_k$ where $U_i=V_i \cap U$ for $i=1, \ldots, k$.

For a subset $U \subset V$ as above, we define a weight 
$$Z^U=\left\{z^U_S: \quad S \in H|U \right\}$$
on $H|U$ as follows. We consider the restriction of weight $Z$ onto hypergraph $H|U$ and define $Z^U$ to be the 
scaling of the restriction to a $k$-stochastic weight.
 We consider the partition function associated with the hypergraph 
 $H|U$, which we denote by $P_U$. We want to estimate $P_U\left(Z^U\right)$.
 
Let $A \in H|U$ be an edge. We consider the complement $U \setminus A$, the corresponding hypergraph 
$H|(U \setminus A)$, weight $Z^{U \setminus A}$ and the partition function $P_{U \setminus A}\left(Z^{U \setminus A}\right)$.

Our goal is to prove that for some $\gamma_1=\gamma_1(k, \alpha)>0$, $\gamma_2=\gamma_2(k, \alpha) >0$
 and $l_0=l_0(k, \alpha)$ we have 
$$\aligned& P_U\left(Z^U\right) \ \geq \ \exp\left\{ -k+1 -{\gamma_1 \over l-1}\right\} \min_{A \in H|U} P_{U\setminus A}\left(Z^{U \setminus A}\right)
 \\
&\qquad \qquad \qquad \text{and} \\
&P_U\left(Z^U\right) \ \leq \ \exp\left\{ -k+1 + {\gamma_2 \over l-1}\right\} \max_{A \in H|U} 
P_{U\setminus A}\left(Z^{U\setminus A}\right) \endaligned \tag5.2.1$$
provided $l \geq l_0$ (recall that $|U|=kl$).

We show that we can choose 
$$\gamma_1 =\alpha^{3(k+1)} (k^2+k)^2 +(k-1)^2, \quad \gamma_2={k^2 \alpha^{k+1} \over 2}  \quad \text{and} \quad 
l_0 = \lceil \alpha^{2(k+1)} k^2 \rceil +1 .$$

Since the restriction of $Z$ onto $H|U$ is $\alpha$-balanced, by Lemma 4.1 the weight $Z^U$ 
 is $\alpha^{k+1}$-balanced. Crude estimates give 
 $$\alpha^{-(k+1)l_0} l_0^{l_0} {kl_0 \choose k}^{1-l_0} \ \leq \ P_U\left(Z^U\right) \ \leq \ \alpha^{(k+1)l_0} l_0^{l_0 -kl_0+k}$$
 if $|U|=kl_0$.
 
 Starting with $U=V$, $l=m$ and $Z^U =Z$, by iterating (5.2.1), we obtain 
$$\split &P(Z) \ \geq \ \alpha^{-(k+1)l_0} l_0^{l_0}
{kl_0 \choose k}^{1-l_0} \exp\left\{ -(k-1)(m-l_0) - \gamma_1 \sum_{j=l_0+1}^m 
{1 \over j-1} \right\} \\
&\qquad \qquad \qquad  \text{and} \\
&P(Z) \ \leq \ \alpha^{(k+1)l_0} l_0^{l_0 -kl_0+k}\exp\left\{ -(k-1)(m-l_0)  +\gamma_2 \sum_{j=l_0+1}^m {1 \over j-1} \right\}.\endsplit$$
In particular, 
$$P(Z)=\exp\left\{ -(k-1)m + O\left(\sum_{j=1}^m {1 \over j}\right) \right\},$$
which completes the proof of the theorem.
 
We proceed to prove (5.2.1), assuming that $l \geq \alpha^{2(k+1)}k^2 +1$. Since weight $Z^U$ is $\alpha^{k+1}$-balanced, we have
 $$z_S^U \ \leq \ \alpha^{k+1}l {kl \choose k}^{-1} \quad \text{for all} \quad S \subset U \qquad \text{when} \quad H={V \choose k} 
 \tag5.2.2$$
 and 
 $$z_S^U \ \leq \ \alpha^{k+1} l^{-k+1} \quad \text{for all} \quad S \subset U \qquad \text{when} \quad H = V_1 \times \ldots \times 
V_k. \tag 5.2.3$$

Let us pick an element $u \in U$. Then there is a recursion 
$$P_U \left(Z^U\right)=\sum \Sb A \in H|U \\ A \ni u \endSb
z^U_A  \cdot P_{U \setminus A} \left(Z^U\right). \tag5.2.4 $$
Here $P_{U\setminus A}\left(Z^U\right)$ is the partition function computed on the restriction of the weight 
$Z^U$ onto the hypergraph $H|(U\setminus A)$, which is not the same as $P_{U \setminus A}\left(Z^{U \setminus A}\right)$,
since the weight $Z^{U\setminus A}$ is obtained from $Z^U$ by restricting it onto $H|(U\setminus A)$ and scaling the restriction
to a $k$-stochastic weight.

Since $Z^U$ is a $k$-stochastic weight on $H|U$, we have 
$$\sum \Sb A \in H|U \\ A \ni u \endSb z_A^U =1 \quad 
\text{and} \quad z_A^U \geq 0 \quad \text{for all} \quad A \in H|U. \tag5.2.5$$
Let 
$$\sigma_A^U =\sum_{S \in H|(U\setminus A) } z^U_S$$
be the sum of the weights in the restriction $Z^U$ onto $H|(U \setminus A)$, that is, the sum of the weights 
$z_S^U$ for the edges $S \subset U$ not intersecting an edge $A \in H|U$. Since 
weight $Z^U$ is $k$-stochastic, we have 
$$\sigma_A^U =l-k + O\left({1 \over l}\right).$$
More precisely, from (5.2.2) we have 
$$l-k \ \leq \ \sigma_A^U \ \leq \ l-k + {k \choose 2} {kl-2 \choose k-2} \alpha^{k+1} l {kl \choose k}^{-1}$$
if $H$ is a complete $k$-uniform hypergraph and from (5.2.3) we have
$$l-k \ \leq \ \sigma_A^U \ \leq \ l -k + {k \choose 2} \alpha^{k+1} l^{-1},$$
if $H$ is a complete $k$-partite hypergraph.
In either case,
$$l- k \ \leq \ \sigma_A^U \ \leq \ l-k +   {\alpha^{k+1} k^2 \over 2(l-1)}. \tag5.2.6$$
Similarly, from (5.2.2) we have 
$$ 1-  k  {kl -2 \choose k-2}\alpha^{k+1} l {kl \choose k}^{-1} \ \leq \ \sum \Sb S \in H|(U \setminus A) \\ S \ni v \endSb z_S^U \ \leq \ 1 \quad \text{for all} \quad v \in U \setminus A$$
if $H$ is a complete $k$-uniform hypergraph
and from (5.2.3) we have 
$$1 -  k \alpha^{k+1} l^{-1} \ \leq \ \sum \Sb S \in H|(U \setminus A) \\ S \ni v \endSb z_S^U \ \leq \ 1 \quad \text{for all} \quad v \in U \setminus A$$
if $H$ is a complete $k$-partite hypergraph.
In either case,
$$1 - { k \alpha^{k+1} \over l-1} \ \leq \ \sum \Sb S \in H|(U \setminus A) \\ S \ni v \endSb z_S^U \ \leq \ 1 \quad \text{for all} \quad v \in U \setminus A. \tag5.2.7$$

Let us define a weight 
$$W^{U \setminus A}=\left\{w^{U \setminus A}_S: \ S \in H|(U\setminus A) \right\}$$ 
 by scaling the restriction of the weight $Z^U$ onto $H|(U\setminus A)$
to the total sum $l-1$, so that 
$$w^{U \setminus A}_S={l-1 \over \sigma_A^U} \ z^U_S \quad \text{for all} 
\quad S \in H|(U\setminus A). $$
Hence
$$P_{U\setminus A}\left(Z^U\right)= 
\left({\sigma_A^U \over l-1}\right)^{l-1} P_{U \setminus A} \left(W^{U \setminus A}\right). \tag5.2.8$$
We have 
$$\left({\sigma_A^U \over l-1}\right)^{l-1}=\exp\left\{ -k+1 + O\left({1 \over l}\right)\right\}.$$
More precisely, from (5.2.6) we have 
$$\exp\left\{-k+1 - {(k-1)^2 \over l-1}\right\} \ \leq\ \left({\sigma_A^U \over l-1}\right)^{l-1} \ \leq \ 
\exp\left\{-k+1 +{k^2 \alpha^{k+1}\over 2(l-1)} \right\}. \tag5.2.9$$

Moreover, from (5.2.6) and (5.2.7) we deduce that
$$1 - { k \alpha^{k+1} \over l-1} \ \leq \ \sum \Sb S \in H|(U \setminus A) \\ S \ni v \endSb w_S^{U\setminus A} \ \leq \ 
1+ {2(k-1) \over l-1} \quad \text{for all} \quad v \in U \setminus A. \tag5.2.10$$

We intend to apply Lemma 4.2 to weight $W^{U \setminus A}$. We observe that weight $W^{U \setminus A}$
is obtained from weight $Z^U$ by restricting it onto the set $U \setminus A$ and then 
scaling to the total sum $l-1$ of components. Therefore, the $k$-stochastic weight on 
$U \setminus A$ obtained from $W^{U \setminus A}$ by scaling is just $Z^{U \setminus A}$, 
the $k$-stochastic weight obtained by restricting the original weight $Z$ onto $U \setminus A$ and
scaling.

We have 
$$P_{U \setminus A} \left( W^{U\setminus A}\right) = \exp\left\{O\left({1 \over l}\right)\right\} P_{U \setminus A} 
\left(Z^{U\setminus A}\right).$$
More precisely, since $W^{U \setminus A}$ is $\alpha^{k+1}$-balanced and (5.2.10) holds, by Lemma 4.2
we conclude that 
$$\aligned &P_{U \setminus A}\left(W^{U \setminus A}\right) \ \geq \ \exp\left\{-{\alpha^{3(k+1)} (k^2+k)^2 \over l-1} \right\} P_{U \setminus A} \left(Z^{U \setminus A} \right) \\ 
&\qquad \qquad \text{and} \\
&P_{U \setminus A}\left(W^{U \setminus A} \right) \ \leq \ P_{U \setminus A} \left(Z^{U \setminus A} \right).
\endaligned \tag5.2.11$$

Combining (5.2.4), (5.2.5), (5.2.8), (5.2.9)  and (5.2.11) we obtain (5.2.1) with 
$$\gamma_1 = \alpha^{3(k+1)} (k^2+k)^2 +(k-1)^2 \quad \text{and} \quad \gamma_2 = {k^2 \alpha^{k+1} \over 2},$$
which completes the proof. 
{\hfill \hfill \hfill} \qed

\subhead (5.3) Proof of Theorem 1.6 \endsubhead 
Let $Z$ be the $k$-stochastic weight obtained from weight $W$ by scaling. 
By Theorem 3.1 we have 
$$P_H(Z) =f_W(Z) P_H(W) =e^{\zeta} P_H(W).$$
Moreover, by Lemma 4.1, weight $Z$ is $\alpha^{k+1}$-balanced 
and the proof follows by Theorem 1.3 applied to $Z$. Furthermore, weights $Z^U$ constructed in Section 5.2, being 
scalings of restrictions of $W$ onto subsets $U \subset V$, are also $\alpha^{k+1}$-balanced, 
and hence we can use the same estimates for $P_H(Z)$ as in Theorem 1.3.
{\hfill \hfill \hfill} \qed

\head Acknowledgment \endhead

The authors are grateful to Jeff Kahn for helpful discussions.

\enddocument
\end